\documentclass[12pt]{amsart}
\usepackage{latexsym,amsmath,amsfonts,amscd,amssymb}
\def\cal{\mathcal}

\newtheorem{theorem}{\bf Theorem}[section]
\newtheorem{lemma}[theorem]{\bf Lemma}
\newtheorem{remark}[theorem]{\bf Remark}

\newtheorem{proposition}[theorem]{\bf Proposition}
\newtheorem{corollary}[theorem]{\bf Corollary}
\newtheorem{definition}[theorem]{\bf Definition}

\newtheorem{compl}[theorem]{\bf Complement}
\newtheorem{claim}[theorem]{\bf Claim}
\newcommand{\be}{\begin{equation}}

\newfont{\bfc}{cmbsy10 scaled 1200}  
\newfont{\dr}{msbm10 scaled \magstep1}  
\newfont{\sdr}{msbm8}  
\newfont{\gl}{eufm10 scaled \magstep1}  

\DeclareFontFamily{OT1}{rsfs}{}
\DeclareFontShape{OT1}{rsfs}{n}{it}{<->rsfs10}{}
\DeclareMathAlphabet{\curly}{OT1}{rsfs}{n}{it}

 \newcommand{\pf}{{\em Proof}. }

 \newcommand{\CC}{{\mathbb C}}

 \newcommand{\HH}{{\mathbb H}}

 \newcommand{\RR}{{\mathbb R}}

 \newcommand{\cD}{{\mathcal D}}
 \newcommand{\cE}{{\mathcal E}}

 \newcommand{\cO}{{\mathcal O}}

 \newcommand{\cU}{{\mathcal U}}

\newcommand{\qu}{/\kern-.7ex/}
\newcommand{\exh}{\to\kern-1.8ex\to}

\newcommand{\End}{\operatorname{End}}
\newcommand{\ra}{\rightarrow}
\newcommand{\lra}{\longrightarrow}

\newcommand{\Hom}{\operatorname{Hom}}
\newcommand{\Ext}{\operatorname{Ext}}

\newcommand{\Aut}{\operatorname{Aut}}

\newcommand{\rk}{\operatorname{rk}}

\begin{document}
\title[Moduli spaces of coherent systems of small slope]{Moduli spaces of 
coherent systems of small slope on algebraic curves}
\thanks{All authors are members of the research group VBAC
(Vector Bundles on Algebraic Curves). Support was received from
a grant from the European Scientific Exchange Programme of the
Royal Society of London and the Consejo
Superior de Investigaciones Cient\'{\i}ficas (15646) and a further grant from 
the Royal Society of London for an International Joint Project (2005/R3).
The first author was partially supported by the National Science
Foundation under grant DMS-0072073. The first, second and fourth authors 
were supported through MEC grant MTM2004-07090-C03-01 (Spain) and the fourth 
also through NSF grant DMS-0111298 (US) during a visit to IAS, Princeton, US. The fifth author was supported by 
the Academia Mexicana de Ciencias during a visit to CIMAT, Guanajuato, 
Mexico.}
\date{10 December 2007}
\author{S.~B.~Bradlow}
  \address{Department of Mathematics\\University of Illinois\\
  Urbana\\IL 61801\\USA}
  \email{bradlow@math.uiuc.edu}
\author{O.~Garc\'{\i}a-Prada}
   \address{Instituto de Ciencias Matem\'aticas CSIC-UAM-UCM-UC3M\\
   Consejo Superior de Investigaciones Cient\'{\i}ficas \\ C/ Serrano, 121\\
   28006 Madrid \\ Spain}
   \email{oscar.garcia-prada@uam.es}
\author{V.~Mercat}
   \address{5 rue Delouvain, 75019 Paris, France}
   \email{mercat@math.jussieu.fr}
\author{V.~Mu\~noz}
   \address{Instituto de Ciencias Matem\'aticas  CSIC-UAM-UCM-UC3M\\
   Consejo Superior de Investigaciones Cient\'{\i}ficas \\ C/ Serrano, 113 bis\\
   28006 Madrid \\ Spain}
\address{Facultad de Matem\'aticas\\ Universidad Complutense de Madrid\\ 
Plaza de Ciencias 3\\ 28040 Madrid\\ Spain}
  \email{vicente.munoz@imaff.cfmac.csic.es}
\author{P.~E.~Newstead}
   \address{Department of Mathematical Sciences \\
   University of Liverpool \\ Peach Street \\
   Liverpool L69 7ZL \\ UK}
   \email{newstead@liverpool.ac.uk}
\begin{abstract} Let $C$ be an algebraic curve of genus $g\ge2$.
A coherent system on $C$ consists of a pair $(E,V)$, where
$E$ is an algebraic vector bundle over $C$ of rank $n$ and
degree $d$ and $V$ is a subspace of dimension $k$ of the
space of sections of $E$. The stability of the coherent
system depends on a parameter $\alpha$.
We study the geometry of the moduli space of coherent
systems for $0<d\leq 2n$. We show that these spaces are irreducible whenever they are non-empty and obtain necessary and sufficient conditions for non-emptiness. 

\end{abstract}

\maketitle
\section{Introduction}\label{section:intro}
Let $C$ be a smooth projective algebraic curve of genus $g\geq 2$. A 
{\em coherent system} on $C$
of {\em type} $(n,d,k)$ is a pair $(E,V)$, where $E$ is a vector bundle on $C$ of rank $n$ and
degree $d$ and $V$ is a linear subspace of the space of sections $H^0(E)$ of 
dimension $k$.
Introduced in \cite{KN}, \cite{RV} and \cite{LeP}, there is a notion of stability for
coherent systems which permits the construction of moduli spaces. This notion depends on a
real parameter $\alpha$, and thus leads to a family of moduli spaces. As described in \cite{BG},
there is a useful relation between these moduli spaces and the Brill-Noether loci in the
moduli spaces of semistable bundles of rank $n$ and degree $d$.

In \cite{BGMN} we began a systematic study of the coherent systems moduli spaces, partly
with a view to applications in higher rank Brill-Noether theory. This study has been continued in \cite{BGMMN}, where we have obtained substantial new information about the geometry and topology of
the moduli spaces for $k\le n$.

For $k>n$, much less is known (but see \cite{T} and, for the case $k=n+1$, \cite{BP} and \cite{BBN}). This case, however,
is of considerable interest, not only in its own right, but also because it is linked with Quot-schemes and with morphisms of $C$ into Grassmannians. The latter connection is a direct generalisation of the classical link between linear systems and morphisms to projective space; coherent systems therefore have an important part to play in understanding the projective geometry of $C$. 

In the present paper, we consider coherent systems with $d\le2n$. Our results are essentially a generalisation and extension of those of \cite{BGN,M1,M2}. More precisely, we show that the 
moduli spaces of $\alpha$-stable coherent systems are irreducible whenever they are non-empty 
and obtain necessary and sufficient conditions for non-emptiness. Even for 
Brill-Noether loci, the irreducibility result is stronger than those 
previously known. The condition for non-emptiness is identical with that for 
Brill-Noether loci except when $C$ is hyperelliptic, $d=2n$ and $k>n$. The 
methods are generally similar to those of \cite{BGN, M1,M2}, except that we make essential use of extensions of coherent systems (rather than simply 
extensions of bundles). In particular, at crucial points in the proof of 
irreducibility and, in the hyperelliptic case, that of non-emptiness, we use the methods of \cite{BGMN} to estimate dimensions of spaces of extensions. 
Some of the estimates are quite delicate and require careful use of these 
methods. We also make use of some of the results of \cite{BGMMN} to handle 
the case $k\le n$, although the present paper can be read independently of 
\cite{BGMMN}.

In order to give full statements of our main results, we need to outline 
some definitions and notations (for more details, see section 
\ref{section:basic}). 
We denote by $G(\alpha;n,d,k)$ the moduli space of $\alpha$-stable coherent 
systems of type $(n,d,k)$, where $\alpha\in\RR$ is a parameter (with 
the necessary conditions $d>0$,  $\alpha>0$ and $\alpha(n-k)<d$ for 
non-emptiness 
of $G(\alpha;n,d,k)$), and by $B(n,d,k)$ the {\em Brill-Noether locus} of 
stable bundles of rank $n$ and degree $d$ with $h^0(E)\ge k$. We write 
$\beta(n,d,k)$ for the ``expected dimension'' of $G(\alpha;n,d,k)$, namely
$$\beta(n,d,k)=n^2(g-1)+1-k(k-d+n(g-1)).$$
Note that the data for defining a coherent system $(E,V)$ can also be 
expressed as an exact sequence
$$0\lra D^*\lra V\otimes\cO\lra E\lra F\oplus T\lra0,$$
where $D$ and $F$ are vector bundles, $T$ is a torsion sheaf and 
$h^0(D^*)=0$. A coherent system $(E,V)$ is said to be {\em generated} if $F$ 
and $T$ are both zero. Finally, we define $U(n,d,k)$ and $U^s(n,d,k)$ by
$$U(n,d,k):=\{(E,V): E \mbox{ is stable and } (E,V) \mbox{ is $\alpha$-stable for } \alpha>0, \alpha(n-k)<d\}$$
and
$$U^s(n,d,k):=\{(E,V): (E,V) \mbox{ is $\alpha$-stable for } \alpha>0, \alpha(n-k)<d\}$$
(see section \ref{section:ne} for further details). Clearly $U(n,d,k)\subset
U^s(n,d,k)$.

We can now state our main results.

\noindent{\bf Theorem \ref{thm:10.3}.}\begin{em}
Suppose that $0<d\leq 2n$ and $\alpha>0$. If $(E,V)\in G(\alpha;n,d,k)$, then $h^0(E^*)=0$. Moreover, if
$G(\alpha;n,d,k)\ne\emptyset$, then it is irreducible, and
 \begin{itemize}
 \item[(a)] if $k<n$, the generic element of $G(\alpha;n,d,k)$ has the form
 $$
 0\rightarrow V\otimes\cO\rightarrow E\rightarrow F\rightarrow0,
 $$
 where $F$ is a vector bundle with $h^0(F^*)=0$;
 \item[(b)] if $k=n$, the generic element of $G(\alpha;n,d,k)$ has the form
 $$
 0\rightarrow V\otimes\cO\rightarrow E\rightarrow T\rightarrow0,
 $$
 where $T$ is a torsion sheaf;
 \item[(c)] if $k>n$, the generic element of $G(\alpha;n,d,k)$ has the form
 $$
 0\rightarrow D^*\rightarrow V\otimes\cO\rightarrow E\rightarrow0,
 $$
 i.e. $(E,V)$ is generated;
 \item[(d)] $\dim G(\alpha;n,d,k)=\beta(n,d,k)$ except when $C$ is 
hyperelliptic and $(n,d,k)=(n,2n,n+1)$ with $n<g-1$.
\end{itemize}\end{em}

\noindent{\bf Theorem \ref{thm:non}.}\begin{em}
Suppose that $C$ is non-hyperelliptic of genus $g\ge3$, $n\ge2$ and $0<d\le2n$. Then
$U(n,d,k)\ne\emptyset$ if and only if either 
$$k\le n+\frac1g(d-n),\ \ (n,d,k)\ne(n,n,n)$$
or 
$$(n,d,k)=(g-1,2g-2,g).$$
In all other cases,
\begin{itemize}
\item $G(\alpha;n,d,k)=\emptyset$ for all $\alpha>0$;
\item $B(n,d,k)=\emptyset$.
\end{itemize}\end{em}

\noindent{\bf Theorem \ref{thm:hyp}.}\begin{em}
Suppose that $C$ is hyperelliptic, $n\ge2$ and $0<d\le2n$. Then
\begin{itemize}
\item[(a)] $U(n,d,k)\ne\emptyset$ if and only if either 
$$0<d<2n,\ \ k\le n+\frac1g(d-n),\ \ (n,d,k)\ne(n,n,n)$$
or $d=2n, k\le n$;
\item[(b)] if $k>n$, then
\begin{itemize}
\item $U(n,2n,k)=\emptyset$;
\item $U^s(n,2n,k)\ne\emptyset$ if and only if either $k\le n+\frac{n}g$ or $k=n+1$ and $2\le n\le g-1$.
\end{itemize}
\end{itemize}
In all other cases,
\begin{itemize}
\item $G(\alpha;n,d,k)=\emptyset$ for all $\alpha>0$;
\item $B(n,d,k)=\emptyset$.
\end{itemize}\end{em}

\noindent The case $n=1$ is omitted from the last two statements since 
the results then need modifying; of course this case is very simple.

The contents of the paper are as follows. In section \ref{section:basic}, 
we give definitions and notations together with some basic facts which we 
shall need. In section \ref{section:le2}, we generalise the results of 
\cite{BGN, M1, M2} to obtain a necessary condition for the existence of 
$\alpha$-stable coherent systems. Section \ref{section:irred} is devoted to 
a proof of irreducibility (Theorem \ref{thm:10.3}). In section 
\ref{section:ne}, we state our results on non-emptiness separately for $C$ 
non-hyperelliptic (Theorem \ref{thm:non}) and for $C$ hyperelliptic 
(Theorem \ref{thm:hyp}); the proofs for $C$ non-hyperelliptic are included.
In the lengthy section \ref{section:hyp} we prove Theorem \ref{thm:hyp}; 
this requires some delicate constructions using the methods of \cite{BGMN}. 
Finally section \ref{section:ex} contains an example with $d>2n$ to show 
that the situation can then be more complicated.

We suppose throughout that $C$ is a smooth projective algebraic curve 
of genus $g\ge2$ defined over the complex numbers. The cases $g=0$ and 
$g=1$ have been investigated in \cite{LN1,LN2,LN3}, where irreducibility has 
been proved with no restriction on the degree, but the non-emptiness 
results for the case $g=0$ are still not complete. We also assume that 
$k\ge1$.

\section{Definitions, notations and basic facts}\label{section:basic}

We refer the reader to \cite{BGMN} for the basic properties of coherent 
systems on algebraic curves. For convenience, we provide here a synopsis 
of the main definitions and facts which we shall need. Recall that the 
{\em slope} $\mu(E)$ of a vector bundle of rank $n$ and degree $d$ is defined
by $\mu(E):=\frac{d}n$.

\begin{definition}\label{def:basic1} \begin{em} Let $(E,V)$ be a coherent system of 
type $(n,d,k)$. For any $\alpha\in\RR$, the {\em $\alpha$-slope}  
$\mu_\alpha(E,V)$ is defined by 
$$\mu_\alpha(E,V):=\frac{d}n+\alpha\frac{k}n.$$ 
A {\em coherent subsystem} of $(E,V)$ is a coherent system $(E',V')$ such 
that $E'$ is a subbundle of $E$ and $V'\subset V\cap H^0(E')$. A {\em 
quotient coherent 
system} of $(E,V)$ is a coherent system $(E'',V'')$ together with a 
homomorphism $(E,V)\ra(E'',V'')$ such that both $E\ra E''$ and $V\ra V''$ 
are surjective.\end{em}
\end{definition}

Note that, with our definition of coherent system, a subsystem possesses 
a corresponding quotient system only if $V'=V\cap H^0(E')$.

\begin{definition}\label{def:basic2}\begin{em}
A coherent system $(E,V)$ is {\em $\alpha$-stable} 
({\em $\alpha$-semistable}) if, for 
every proper coherent subsystem $(E',V')$,
$$\mu_\alpha(E',V')<(\le)\mu_\alpha(E,V).$$
\end{em}\end{definition}
There exists a moduli space $G(\alpha;n,d,k)$ of $\alpha$-stable coherent 
systems of type $(n,d,k)$; necessary conditions for non-emptiness are 
$$d>0,\ \ \alpha>0,\ \ (n-k)d<\alpha.$$

\begin{definition}\label{def:basic3}\begin{em}
A {\em critical value} for coherent systems of type $(n,d,k)$ is a value 
of $\alpha>0$ for which there exists a coherent system $(E,V)$ of type  
$(n,d,k)$ and a coherent subsystem $(E',V')$ of $(E,V)$ of type 
$(n',d',k')$ such that $\frac{k'}{n'}\ne\frac{k}n$ but 
$\mu_\alpha(E',V')=\mu_\alpha(E,V)$. We also regard $\alpha=0$ as a critical value. 
\end{em}\end{definition}

It is known \cite[Propositions 4.2, 4.6]{BGMN} that, for any $(n,d,k)$, 
there are finitely many critical values
$$0=\alpha_0<\alpha_1<\ldots<\alpha_L<\left\{\begin{array}{ll}
\frac{d}{n-k}&\mbox{ if } k<n\\\infty&\mbox{ if } k\ge n.\end{array}\right.$$
Moreover, if $k<n$ and $\alpha\ge\frac{d}{n-k}$, then 
$G(\alpha;n,d,k)=\emptyset$. For $\alpha,\alpha'\in(\alpha_i,\alpha_{i+1})$, 
we have $G(\alpha;n,d,k)=G(\alpha';n,d,k)$ and we denote this moduli space 
by $G_i:=G_i(n,d,k)$. We shall be particularly concerned with the moduli 
spaces $G_0$ (``small'' $\alpha$) and $G_L$ (``large'' $\alpha$). If 
$(E,V)\in G_0$, we say also that $(E,V)$ is {\em $0^+$-stable} (with similar 
definitions for {\em $\alpha^\pm$-stable}).

We denote by $M(n,d)$ the moduli space of stable bundles of rank $n$ and 
degree $d$, and by $B(n,d,k)$ the {\em Brill-Noether locus}
$$B(n,d,k):=\{E\in M(n,d):h^0(E)\ge k\}.$$
We have, for any coherent system $(E,V)$, \cite[Proposition 2.5]{BGMN}
\begin{itemize}
\item $(E,V)\in G_0 \Longrightarrow E\mbox{ semistable}$;
\item $E\mbox{ stable } \Longrightarrow (E,V)\in G_0$.
\end{itemize} 

The moduli space $G(\alpha;n,d,k)$ has the property that every irreducible 
component has dimension greater than or equal to the {\em Brill-Noether 
number} 
\begin{equation}\label{eqn:basic1}
\beta(n,d,k):=n^2(g-1)+1-k(k-d+n(g-1)).
\end{equation}
This number is the ``expected dimension'' of $G(\alpha;n,d,k)$ in a stronger 
sense. For this, we define, 
for any coherent system $(E,V)$, the {\em Petri map} of $(E,V)$ as the map
$$V\otimes H^0(E^*\otimes K)\lra H^0(E\otimes E^*\otimes K)$$
given by multiplication of sections. This map governs the infinitesimal 
behaviour of the moduli space in the following sense \cite[Proposition 3.10]{BGMN}:
\begin{itemize}
\item Let $(E,V)$ be an $\alpha$-stable coherent system of type $(n,d,k)$. 
Then $G(\alpha;n,d,k)$ is smooth of dimension $\beta(n,d,k)$ at the point 
corresponding to $(E,V)$ if and only if the Petri map of $(E,V)$ is 
injective.
\end{itemize}.

\begin{definition}\label{def:basic4}\begin{em}
The coherent system $(E,V)$ is {\em generated} if the evaluation map 
$V\otimes\cO\ra E$ is surjective. The bundle $E$ is {\em generated} if 
$(E,H^0(E))$ is generated.\end{em}
\end{definition}

We shall make no explicit use of the flip loci $G_i^\pm$ of \cite[section 6]
{BGMN}, so shall not describe them here. However we make extensive use 
of extensions
\begin{equation}\label{eqn:basic2}
0\lra(E_1,V_1)\lra(E,V)\lra(E_2,V_2)\lra0,
\end{equation}
where $(E_1,V_1)$, $(E_2,V_2)$ are coherent systems of types $(n_1,d_1,k_1)$,
$(n_2,d_2,k_2)$ respectively. Here we use the notations and results of 
\cite[section 3]{BGMN}. The extensions (\ref{eqn:basic2}) are classified in 
the usual way by a group 
$$\Ext^1((E_2,V_2),(E_1,V_1)).$$
For dimensional reasons $\Ext^q((E_2,V_2),(E_1,V_1))=0$ for $q\ge3$, so 
we have \cite[equation (8)]{BGMN}
\begin{equation}\label{eqn:basic3}
\dim\Ext^1((E_2,V_2),(E_1,V_1))=C_{21}+\dim\HH^0_{21}+\dim\HH^2_{21},
\end{equation}
where
\begin{equation}\label{eqn:basic4}
C_{21}:=n_1(n_2-k_2)(g-1)+(k_2-n_2)d_1+d_2n_1-k_1k_2
\end{equation}
and
$$\HH^0_{21}=\Hom((E_2,V_2),(E_1,V_1)),\ \ 
\HH^2_{21}=\Ext^2((E_2,V_2),(E_1,V_1)).$$
The main purpose of introducing the number $C_{21}$ is that frequently, 
although not always, $\HH^0_{21}$ and $\HH^2_{21}$ are both zero and 
$\dim\Ext^1((E_2,V_2),(E_1,V_1))$ is then given by the purely numerical 
formula (\ref{eqn:basic4}). Of course, we can define $\HH^0_{12}$, 
$\HH^2_{12}$ and $C_{12}$ by interchanging the indices, and in 
particular
\begin{equation}\label{eqn:basic5}
C_{12}=n_2(n_1-k_1)(g-1)+(k_1-n_1)d_2+d_1n_2-k_1k_2.
\end{equation}
Note \cite[Corollary 3.7]{BGMN} that, with the notation of 
(\ref{eqn:basic2}), 
\begin{equation}\label{eqn:basic6}
\beta(n,d,k)=\beta(n_1,d_1,k_1)+\beta(n_2,d_2,k_2)+C_{12}+C_{21}-1.
\end{equation}
Note further \cite[equation (11)]{BGMN} that, if $N_2$ is the kernel of the 
evaluation map $V_2\otimes\cO\ra E_2$, then
\begin{equation}\label{eqn:basic7}
\HH^2_{21}\cong H^0(E_1^*\otimes N_2\otimes K)^*.
\end{equation}

Putting $(E_1,V_1)=(E_2,V_2)$ in (\ref{eqn:basic3}) and using 
(\ref{eqn:basic1}) and (\ref{eqn:basic4}), we get
$$\dim\Ext^1((E,V),(E,V))=\beta(n,d,k)+\dim\End(E,V)+
\dim\Ext^2((E,V),(E,V))-1.$$
Now, when $\Ext^2((E,V),(E,V))=0$, there is no obstruction to the 
construction of a local deformation space for $(E,V)$ and this local 
deformation space has dimension
\begin{equation}\label{eqn:basic8}
\dim\Ext^1((E,V),(E,V))=\beta(n,d,k)+\dim\End(E,V)-1
\end{equation}
(see the proof of \cite[Th\'eor\`eme 3.12]{He}). Note further

\begin{proposition}\label{prop:extra}
$\Ext^2((E,V),(E,V))=0$ if and only if the Petri map of $(E,V)$ is injective.
\end{proposition}
\pf By (\ref{eqn:basic7}), $\Ext^2((E,V),(E,V))\cong H^0(E^*\otimes N\otimes K)^*$, 
where $N$ is the kernel of the evaluation map $V\otimes\cO\ra E$. It follows from the evaluation sequence $0\ra N\ra V\otimes\cO\ra E$ that  $H^0(E^*\otimes N\otimes K)$ is the kernel of the Petri map, giving the result.\qed

We need one further important fact about the extensions (\ref{eqn:basic2}).
\begin{proposition}\label{prop:basic}
If (\ref{eqn:basic2}) is non-trivial and $\alpha_i$ is a critical value 
such that $(E_1,V_1)$ and $(E_2,V_2)$ are both $\alpha_i$-stable with 
$\mu_{\alpha_i}(E_1,V_1)=\mu_{\alpha_i}(E_2,V_2)$ and 
$\mu_{\alpha_i^-}(E_1,V_1)<\mu_{\alpha_i^-}(E_2,V_2)$, then $(E,V)$ is 
$\alpha_i^-$-stable.
\end{proposition}

Since this is not explicitly stated in either \cite{BGMN} or \cite{BGMMN} 
(although it is used in \cite{BGMMN}), we give a proof.

\pf Suppose that $(E',V')$ is a coherent subsystem of $(E,V)$ contradicting 
$\alpha_i^-$-stability. Then $(E',V')$ also contradicts $\alpha_i$-stability 
of $(E,V)$. It follows that either $(E',V')=(E_1,V_1)$ or $(E',V')$ maps 
isomorphically to $(E_2,V_2)$. In the first case, $\alpha_i^-$-stability of 
$(E,V)$ is not contradicted, while in the second (\ref{eqn:basic2}) is 
trivial.\qed

\section{Coherent systems for $d\le2n$}\label{section:le2}

Our first object in this section is
to obtain a necessary condition for the existence of $\alpha$-stable coherent
systems for $d\le2n$; it turns out that the condition is almost identical with
that for stable bundles (see \cite{M1, M2}).

We start with the case $d<2n$, when the results of \cite{M1} carry
over quite easily to give a
necessary condition for $\alpha$-semistability.

\begin{lemma}\label{lemma:prelim}
Suppose that $(E,V)$ is an $\alpha$-semistable coherent system for
some $\alpha>0$ and that $0<d<2n$. Then
\begin{equation}\label{eq:kbound}
k\le n+\frac1g(d-n).
\end{equation}\end{lemma}

\pf If $E$ is semistable, the result holds by \cite[Chapitre 2,
Th\'eor\`eme A.1]{M1}.

If $E$ is not semistable, then $E$ has a stable quotient $G$ with
$\mu(G)<\mu(E)<2$. Again by \cite[Chapitre 2, Th\'eor\`eme
A.1]{M1}, we have
$$h^0(G)\le n_G+\frac1g(d_G-n_G),
$$
where $n_G$ and $d_G$ denote the rank and degree of $G$. Let $W$
denote the image of $V$ in $H^0(G)$. Then, if $k>n+\frac1g(d-n)$,
we have
$$
\frac{\dim
W}{n_G}\le1+\frac1g\left(\frac{d_G}{n_G}-1\right)<1+\frac1g\left(\frac{d}n-1\right)<\frac{k}n.
$$
It follows that the quotient coherent system $(G,W)$ contradicts
the $\alpha$-semistability of $(E,V)$ for any $\alpha>0$. \qed

This lemma has the following interesting consequence, which has
relevance for coherent systems with $k>n$ in general.

\begin{corollary}\label{cor:k>n}
Suppose that $k>n$ and that there exists an $\alpha$-semistable
coherent system $(E,V)$ for some $\alpha>0$. Then
$$
d\ge\min\left\{2n,n+g(k-n)\right\}.
$$
\end{corollary}

\pf If $k>n$ and $0<d<2n$, then the lemma implies that $d-n\ge
g(k-n)$. For $d=0$, $(E,V)$ cannot be $\alpha$-stable. The associated graded object must be a sum of
coherent systems of types $(1,0,1)$ or $(1,0,0)$ and hence $k\le
n$. \qed

In order to cover the case $d=2n$, we shall make use of the dual
span construction, which we briefly recall. Let $(F,W)$ be a
coherent system. Slightly modifying the notations of \cite[section
5.4]{BGMN}, we define a coherent system
$$
D(F,W)=(D_W(F),W'),
$$
where $D_W(F)$ is defined by the exact sequence
$$
0\longrightarrow D_W(F)^*\lra W\otimes {\cal O}\longrightarrow F
$$
and $W'$ is the image of $W^*$ in $H^0(D_W(F))$. In particular we
write $D(F)$ for $D_{H^0(F)}(F)$. If $W\otimes\cO\ra F$ is
surjective, then the linear map $W^*\ra W'$ is induced from the
dual exact sequence \be\label{eq:dualf} 0\lra F^*\lra
W^*\otimes\cO\lra D_W(F)\lra0.
\end{equation}
Moreover, if $h^0(F^*)=0$, then $W^*$ maps isomorphically to $W'$.

For the canonical line bundle $K$, we obtain a bundle $D(K)$ of
rank $g-1$ and degree $2g-2$. Taking $W=H^0(K)$, (\ref{eq:dualf})
becomes \be\label{eq:dualk} 0\lra K^*\lra H^0(K)^*\otimes\cO\lra
D(K)\lra0.
\end{equation}
It is known \cite[Corollary 3.5 and Remark 3.6(i)]{PR} that, if $C$ is not hyperelliptic, then
$D(K)$ is stable and $h^0(D(K))=g$, while, if $C$ is
hyperelliptic, then $D(K)\cong L^{\oplus(g-1)}$, where $L$ is the
hyperelliptic line bundle. In both cases, we obtain new
$\alpha$-stable coherent systems with $d=2n$; to describe them, we
use the following general lemma.

\begin{lemma}\label{lemma:all}
Let $(E,V)$ be a generated coherent system of type $(n,d,n+1)$ such that $E$
is a semistable bundle. Then $(E,V)$ is
$\alpha$-stable for all $\alpha>0$.
\end{lemma}

\pf Let $(F,W)$ be a coherent subsystem of $(E,V)$ with $0<\rk
F<n$. Since $E$ is semistable, $\mu(F)\le\mu(E)$. To show that
$\mu_\alpha(F,W)<\mu_\alpha(E,V)$, it is therefore sufficient to
show $\dim W\le\rk F$. Suppose $\dim W>\rk F$. Then the image of
$V$ in $H^0(E/F)$ has dimension $\le\rk(E/F)$. Now $d>0$, hence
$\deg(E/F)>0$. It follows that the image of $V$ does not generate
$E/F$, which contradicts the hypothesis that $V$ generates $E$.
\qed

\begin{corollary}\label{cor:nhyp}
If $C$ is not hyperelliptic, then $D(K,H^0(K))$ is $\alpha$-stable
of type $(g-1,2g-2,g)$ for all $\alpha>0$.
\end {corollary}

\pf This follows immediately from the lemma and (\ref{eq:dualk}).
\qed

\begin{corollary}\label{cor:hyp}
Suppose $C$ is hyperelliptic and $a$ is an integer, $1\le a\le
g-1$. Let $L$ be the hyperelliptic line bundle and
$W$ a subspace of $H^0(L^{\oplus a})$ of dimension
$a+1$ which generates $L^{\oplus a}$. Then
\begin{itemize}
\item the coherent system $(L^{\oplus a},W)$ is $\alpha$-stable
of type $(a,2a,a+1)$ for all $\alpha>0$; \item $(L^{\oplus a},W)\cong
D(L^a,H^0(L^a))$.
\end{itemize}
\end{corollary}

\pf The first statement follows at once from the lemma. For the
second statement, note that we have an exact sequence
$$
0\lra M\lra W\otimes\cO\lra L^{\oplus a}\lra0,
$$
where $M$ is a line bundle. But then $M^*\cong\det L^{\oplus
a}\cong L^a$. We therefore have a dual exact sequence
$$
0\lra (L^{\oplus a})^*\lra W^*\otimes\cO\lra L^a\lra0.
$$
Since $h^0(L^a)=a+1$, this must be the defining exact sequence for
$D(L^a,H^0(L^a))$, which completes the proof. \qed

\begin{remark}\label{rmk8}
{\em In Corollary \ref{cor:hyp}, for any subspace $W$ of dimension $a+1$ which generates $L^{\oplus a}$, the isomorphism class of $(L^{\oplus a},W)$ is the same.}
\end{remark}

\begin{lemma}\label{lemma:csmorphism}
Let $(E,V)$ be a coherent system and $F$ a vector bundle. Suppose that $F$ is
generated and that $h^0(F^*)=0$. Then
$$
\Hom(D(F,H^0(F)),(E,V))
$$
is isomorphic to the kernel of the homomorphism
$$
\Psi:H^0(F)\otimes V\rightarrow H^0(F\otimes E)
$$
given by multiplication of sections.
\end{lemma}
\pf We have an exact sequence of coherent systems
\be\label{eq:csseq} 0\longrightarrow (F^*,0)\longrightarrow
(H^0(F)^*\otimes \cO,H^0(F)^*)\longrightarrow
D(F,H^0(F))\longrightarrow 0\ .
\end{equation}
Taking $\Hom((\ref{eq:csseq}),(E,V))$, we get an exact sequence
\begin{eqnarray*}
0\lra\Hom(D(F,H^0(F)),(E,V))&\lra&\Hom((H^0(F)^*\otimes\cO,H^0(F)^*),(E,V))\\
&\stackrel{\psi}{\lra}&\Hom((F^*,0),(E,V))\lra\ldots
\end{eqnarray*}
Now $\psi$ can be identified with the natural linear map
$$
\Hom((H^0(F)^*,V)\lra\Hom(F^*,E)
$$
and this in turn can be identified with $\Psi$. \qed

\begin{corollary}\label{cor:stable}
Let $(E,V)$ be a coherent system of type $(n,d,k)$ with
$h^0(E^*)=0$ and let
$$
m=\dim\Hom(D(K,H^0(K)),(E,V)).
$$
Then
$$
k\le n+\frac1g(d-n+m).
$$
\end{corollary}

\pf Apply Lemma \ref{lemma:csmorphism} with $F=K$.
The condition $h^0(E^*)=0$ implies by Serre duality that
$h^0(K\otimes E)=d+n(g-1)$. \qed

\begin{remark}\label{rmk7}
{\em If $E$ is a semistable
bundle with $0<d<2n$, then $\Hom(D(K),E)=0$ since $D(K)$ is semistable of slope $2$. So $m=0$ and the corollary reduces to
\cite[Chapitre 2, Th\'eor\`eme A.1]{M1}.}
\end{remark}

We come now to the main result of this section. Although we have
already proved it in the case $d<2n$, for completeness we state it
for the whole range $d\le 2n$.

\begin{proposition}\label{prop:bound}
Let $(E,V)$ be an $\alpha$-stable coherent system of type $(n,d,k)$
with $0<d\leq2n$. Then
$$
k\leq n+\frac{1}{g}(d-n)
$$
except when $d=2n$ and one of the following holds:
\begin{itemize}
\item $C$ is not hyperelliptic and $(E,V)\cong D(K,H^0(K))$;
\item $C$ is hyperelliptic and $(E,V)\cong(L^{\oplus a},W)$,
where $L$ is the hyperelliptic line bundle, $a\le g-1$ and $W$
is a subspace of $H^0(L^{\oplus a})$ of dimension $a+1$
which generates $L^{\oplus a}$.
\end{itemize}
\end{proposition}
\pf For $d<2n$, this follows at once from Lemma
\ref{lemma:prelim}. So we can suppose $d=2n$. If $E$ is stable,
the proposition follows from the results of \cite{M2}. If $E$ is
not semistable, the proof of Lemma \ref{lemma:prelim} still works.

It remains to consider the case where $E$ is strictly semistable
with $d=2n$. We can certainly suppose that \be\label{eq:knn}
k>n+\frac{1}{g}(d-n)=n\left(1+\frac1g\right).
\end{equation}
By Corollary \ref{cor:stable} this implies that there exists a
non-zero homomorphism \be\label{eq:hom} D(K,H^0(K))\lra(E,V).
\end{equation}

Suppose first that $C$ is not hyperelliptic. Then $D(K)$ is
stable; since $E$ is strictly semistable, the homomorphism
(\ref{eq:hom}) must be injective and indeed
$$
\dim\Hom(D(K,H^0(K)),(E,V))\le\dim\Hom(D(K),E)\le\frac{n}{g-1}.
$$
Corollary \ref{cor:stable} implies that
$$
k\le n\left(1+\frac1g\right)+\frac{n}{g(g-1)}=\frac{ng}{g-1}.
$$
But now
$$
\mu_{\alpha}(D(K,H^0(K)))=2+\alpha\frac{g}{g-1}\ge2+\alpha\frac{k}n,
$$
which contradicts the $\alpha$-stability of $(E,V)$ unless $(E,V)\cong D(K,H^0(K))$.\\

If $C$ is hyperelliptic, we have $D(K)\cong L^{\oplus(g-1)}$ and
$h^0(L^*\otimes E)\le n$, so Corollary \ref{cor:stable} gives
$$
 k\leq n\left(1+\frac1g\right)+\frac1gn(g-1)= 2n.
$$
By (\ref{eq:knn}), we deduce that there exists an integer $a$,
$1\le a\le g-1$ such that
$$
n\left(1+\frac{1}{a}\right)\geq k>n\left(1+\frac{1}{a+1}\right).
$$
By Clifford's Theorem (see \cite[Theorem 2.1]{BGN}), we have $h^0(E\otimes L^a)\leq (a+2)n$; so
$$
k\cdot h^0(L^a)=k(a+1)>n(a+2)\ge h^0(E\otimes L^a).
$$
Hence, by Lemma \ref{lemma:csmorphism}, there exists a non-zero
homomorphism of coherent systems
$$
D(L^a,H^0(L^a))\lra(E,V).
$$
Now
$$
\mu_\alpha(D(L^a,H^0(L^a)))=2+\alpha\frac{a+1}a\ge2+\alpha\frac{k}n.
$$
By Corollary \ref{cor:hyp}, this contradicts the
$\alpha$-stability of $(E,V)$ unless
$$
(E,V)\cong D(L^a,H^0(L^a))\cong(L^{\oplus a},W),
$$
where $W$ is any subspace of $H^0(L^{\oplus a})$ of dimension
$a+1$ which generates $L^{\oplus a}$. \qed

\section{Irreducibility of the moduli space for $d\le2n$}\label{section:irred}

In this section we
prove that the moduli space $G(\alpha;n,d,k)$ is
irreducible for $0<d\leq2n$. We start with two lemmas.

\begin{lemma} \label{lem:h1}
 Suppose that $(E,V)$ is a coherent system of
 type $(n,d,k)$ and consider the exact sequence
 \begin{equation}\label{eqn:lemh1}
 0\rightarrow D^*\rightarrow V\otimes {\cal O}\rightarrow
 E\rightarrow F\oplus T\rightarrow 0,
 \end{equation}
 where $D=D_V(E)$, $T$ is a torsion sheaf and $F$ is a vector bundle. 
Suppose further that $\Hom(D(K,H^0(K)),(E,V))=0$. Then
\begin{itemize}
 \item[(a)] $h^1(D)=0$;
\item[(b)] if $F=0$, the Petri map of $(E,V)$ is injective.
\end{itemize}
\end{lemma}

\pf (a) Suppose that $h^1(D)\neq 0$. Then there is a non-zero
homomorphism $D\rightarrow K$. Since $V^*$ generates $D$,
the map $V^*\rightarrow H^0(K)$ is non-zero. Dualising we obtain a diagram

$$
\begin{array}{ccccccccc}
 0&\rightarrow&K^*&\rightarrow&H^0(K)^*\otimes{\cal O} &\rightarrow& D(K)&\rightarrow&0\\
 &&\downarrow&&\downarrow& &\downarrow&& \\
 0&\rightarrow&D^*&\rightarrow&V\otimes {\cal O}
 &\rightarrow&E && \\
 \end{array}
 $$
and hence a non-zero homomorphism $D(K,H^0(K))\rightarrow (E,V)$, contradicting 
the hypothesis.

(b) Tensor the sequence $V \otimes {\mathcal{O}} \rightarrow
E\rightarrow T\rightarrow 0$ with $D$ and apply cohomology to
get an exact sequence
$$
H^1(D\otimes V)\rightarrow H^1(D\otimes E)\rightarrow0.
$$
By (a), $h^1(D)=0$, which implies that $h^1(D\otimes V)=0$, hence $h^1(D\otimes E)=0$. Now dualise the sequence
(\ref{eqn:lemh1}), with $F=0$, to get
 \begin{equation}\label{eqn:h2}
 0\rightarrow E^*\rightarrow V^*\otimes{\cal O}\rightarrow
 D\oplus T\rightarrow 0.
 \end{equation}
Tensor this exact sequence with $E$ to obtain the surjectivity of
the map $H^1(E^*\otimes E) \rightarrow V^* \otimes H^1(E)$.
This map is dual to the Petri map at $(E,V)$, which is
therefore injective. \qed

\begin{lemma} \label{lem:h2}
 Let $\alpha>0$, $d\leq 2n$ and $k\leq n(1+\frac{1}{g})$.
 Let $(E,V)$ be an $\alpha$-semistable coherent system of
 type $(n,d,k)$. Then
 $$\Hom(D(K,H^0(K)),(E,V))=0.$$
\end{lemma}

\pf We have
$$\mu_{\alpha}(D(K,H^0(K)))=2+\alpha\frac{g}{g-1}>2+\alpha
\left(1+\frac1g\right)\ge\mu_\alpha(E,V).
$$
Since $D(K,H^0(K))$ and $(E,V)$ are both $\alpha$-semistable,
this implies that
$$\Hom(D(K,H^0(K)),(E,V))=0.$$\qed

\

Now let $(E,V)$ be a coherent system and let $E'$ be the (subsheaf) image of the evaluation map $V\otimes\cO\rightarrow E$. If $h^0(E'^*)\ne0$, there exists a non-zero homomorphism $E'\ra\cO$. Since $V\otimes\cO\ra E'$ is surjective, this 
induces a non-zero homomorphism $V\otimes\cO\ra\cO$, which necessarily splits; so 
$E'$ admits $\cO$ as a direct summand. It follows by induction that we can write $E'=\cO^{k_2}\oplus G$ 
where $h^0(G^*)=0$. We have a diagram (extending the sequence (\ref{eqn:lemh1}))

\begin{equation}\label{array}
\begin{array}{ccccccccccc}
 &&0&&0&&0&&&&\\
 &&\downarrow&&\downarrow& &\downarrow&&&& \\
 0&\rightarrow&D^*&\rightarrow&V_1\otimes\cO&\rightarrow& E_1 &\rightarrow&T_1&\rightarrow 0\\
 &&\downarrow&&\downarrow& &\downarrow&&\downarrow&& \\
 0&\rightarrow&D^*&\rightarrow& V\otimes{\cal O}&\rightarrow&E
 &\rightarrow& F\oplus T &\rightarrow& 0\\
 &&\downarrow&&\downarrow& &\downarrow&&\downarrow&& \\
 &&0&\rightarrow&\cO^{k_2}&\rightarrow&E_2 &\rightarrow&F\oplus T_2 &\rightarrow& 0\\
 &&&&\downarrow& &\downarrow&&\downarrow&& \\
 &&&&0& &0&&0&& \\
 \end{array}
 \end{equation}
with exact rows and columns. Here $F$ is a  vector bundle and $T$, $T_1$, $T_2$ are torsion sheaves; 
moreover $G$ is the image of $V_1\otimes{\mathcal O}$ in $E$ and $E_1$ is the saturation of $G$ in $E$. 
Writing $V_2:=H^0(\cO^{k_2})\subset H^0(E_2)$, we can interpret (\ref{array}) 
as an exact sequence of coherent systems
\begin{equation}\label{eqn:seq}
0\rightarrow(E_1,V_1)\rightarrow(E,V)\rightarrow(E_2,V_2)\rightarrow0.
\end{equation}

\begin{lemma}\label{lemma:<beta}
Let $0<d\le2n$ and $k\le n+\frac1g(d-n)$. Suppose that $(E_1,V_1)$, $(E_2,V_2)$ are of fixed types $(n_1,d_1,k_1)$, $(n_2,d_2,k_2)$ with $E_1\ne0$, $E_2\ne0$, $h^0(D^*)=0$. Suppose further that $h^0(G^*)=0$, where $G$ is the (sheaf-theoretic) image of $V_1\otimes\cO$ in $E_1$. Then the diagrams (\ref{array}) in which $(E,V)$ is $\alpha$-stable for some $\alpha>0$ depend on fewer than $\beta(n,d,k)$ parameters.
\end{lemma}
\pf
By Lemmas \ref{lem:h1}(a) and \ref{lem:h2}, $h^1(D)=0$. Hence, from the 
cohomology sequence associated to the top row of (\ref{array}), 
$h^1(D\otimes E_1)=0$; thus, by (\ref{eqn:basic7}) and Serre 
duality,
$$\Ext^2((E_1,V_1),(E_1,V_1))\cong H^0(E_1^*\otimes D^*\otimes K)^*=0.$$ 
It follows by (\ref{eqn:basic8}) that the local deformation space of $(E_1,V_1)$ has dimension
$$x_1:=\beta(n_1,d_1,k_1)+\dim\End(E_1,V_1)-1.$$
On the other hand, the Petri map of $(E_2,V_2)$ is clearly injective, so 
by (\ref{eqn:basic8}) and Proposition \ref{prop:extra}, the local deformation space of $(E_2,V_2)$ 
has dimension
$$x_2:=\beta(n_2,d_2,k_2)+\dim\End(E_2,V_2)-1.$$

We need to consider only those extensions (\ref{eqn:seq}) for which 
$(E,V)$ is $\alpha$-stable for some $\alpha$. For fixed $(E_1,V_1)$, 
$(E_2,V_2)$, the group
$$\Aut(E_1,V_1)\times \Aut(E_2,V_2)/\{(\lambda,\lambda^{-1}):\lambda\in\CC^*\}$$
acts freely on these extensions. Hence, in (\ref{eqn:seq}), $(E,V)$ depends on at most
\begin{eqnarray}\label{eqn:atmost}
x_1+x_2&+&\dim\Ext^1((E_2,V_2),(E_1,V_1))-(\dim\Aut(E_1,V_1)+\dim\Aut(E_2,V_2)-1)\nonumber\\
&=&\beta(n_1,d_1,k_1)+\beta(n_2,d_2,k_2)+\dim\Ext^1((E_2,V_2),(E_1,V_1))-1
\end{eqnarray}
parameters. Now, by (\ref{eqn:basic3}), we have
$$\dim\Ext^1((E_2,V_2),(E_1,V_1))=C_{21}+\dim\HH^0_{21}+\dim\HH^2_{21}.$$
Here $\HH^0_{21}=\Hom((E_2,V_2),(E_1,V_1))=0$, since the existence of a non-zero homomorphism would imply that $(E,V)$ is not simple, in contradiction to \cite[Proposition 2.2(ii)]{BGMN}. Moreover, by (\ref{eqn:basic7}),
$$\HH^2_{21}\cong H^0(E_1^*\otimes N_2\otimes K)^*,$$
where $N_2$ is the kernel of the evaluation map $V_2\otimes\cO\rightarrow E_2$, which is clearly $0$. So (\ref{eqn:atmost}) becomes
\begin{equation}\label{eqn:atmost1}
\beta(n_1,d_1,k_1)+\beta(n_2,d_2,k_2)+C_{21}-1.
\end{equation}
So, to prove that the number given by (\ref{eqn:atmost1}) is less than 
$\beta(n,d,k)$, it is enough by (\ref{eqn:basic6}) to prove that 
$C_{12}\ge1$. Now, by (\ref{eqn:basic5}),
 \begin{eqnarray}\label{eqn:lambd}
  C_{12} &=& n_2(n_1-k_1)(g-1)+(k_1-n_1)d_2+d_1n_2-k_1k_2
  \nonumber\\
  &=&(d_1-n_1+(n_1-k_1)g)n_2 +k_1 (n_2-k_2)+ d_2(k_1-n_1).
 \end{eqnarray}
We can now check that the third term in (\ref{eqn:lambd}) is positive and the other two are non-negative.
\begin{itemize}
\item Since $h^0(G^*)=0$ and $G\ne0$, we have $k_1>n_1$ and $d_1>0$; also $k_2\le n_2$, hence $\frac{k_1}{n_1}>\frac{k_2}{n_2}$. Now $\alpha$-stability of $(E,V)$ implies that $\frac{d_1}{n_1}<\frac{d_2}{n_2}$, hence $d_2>\frac{d_1n_2}{n_1}>0$. So $d_2(k_1-n_1)>0$.
\item Since $k_2\le n_2$, $k_1(n_2-k_2)\ge0$.
\item Since Lemma \ref{lem:h2} applies to $(E,V)$, it follows from (\ref{eqn:seq}) that $$\Hom(D(K,H^0(K)),(E_1,V_1))=0.$$ Since $E_1/G$ is a torsion sheaf, $E_1^*$ is a subsheaf of $G^*$. Hence $h^0(E_1^*)=0$ 
and it follows from Corollary \ref{cor:stable} that $k_1\le n_1+\frac1g(d_1-n_1)$, so $d_1-n_1+(n_1-k_1)g\ge0$.
\end{itemize}
This completes the proof of the lemma.
\qed

\begin{theorem} \label{thm:10.3}
Suppose that $0<d\leq 2n$ and $\alpha>0$. If $(E,V)\in G(\alpha;n,d,k)$, then $h^0(E^*)=0$. Moreover, if
$G(\alpha;n,d,k)\ne\emptyset$, then it is irreducible, and
 \begin{itemize}
 \item[(a)] if $k<n$, the generic element of $G(\alpha;n,d,k)$ has the form
 \begin{equation}\label{eqn:k<n}
 0\rightarrow V\otimes\cO\rightarrow E\rightarrow F\rightarrow0,
 \end{equation}
 where $F$ is a vector bundle with $h^0(F^*)=0$;
 \item[(b)] if $k=n$, the generic element of $G(\alpha;n,d,k)$ has the form
 \begin{equation}\label{eqn:k=n}
 0\rightarrow V\otimes\cO\rightarrow E\rightarrow T\rightarrow0,
 \end{equation}
 where $T$ is a torsion sheaf;
 \item[(c)] if $k>n$, the generic element of $G(\alpha;n,d,k)$ has the form
 \begin{equation}\label{eqn:k>n}
 0\rightarrow D^*\rightarrow V\otimes\cO\rightarrow E\rightarrow0,
 \end{equation}
 i.e. $(E,V)$ is generated;
 \item[(d)] $\dim G(\alpha;n,d,k)=\beta(n,d,k)$ except when $C$ is 
hyperelliptic and $(n,d,k)=(n,2n,n+1)$ with $n<g-1$.
\end{itemize}
\end{theorem}

\pf For the fact that $h^0(E^*)=0$, see \cite[Lemma 2.9]{BGMMN}.

Let $Z$ be a component of $G(\alpha;n,d,k)$. Note that $\dim Z\ge\beta(n,d,k)$.

(a) Suppose $k<n$. By Lemma \ref{lemma:<beta}, the generic element of $Z$ must have $(E_1,V_1)=0$ in (\ref{eqn:seq}), so we have
$$0\rightarrow V\otimes\cO\rightarrow E\rightarrow F\oplus T\rightarrow0,$$
i.e. $(E,V)$ is injective in the sense of \cite[Definition 2.1]{BGMMN}. The result now follows from \cite[Theorem 3.3(iii) and its proof]{BGMMN}.

(b) Suppose $k=n$. By Lemma \ref{lemma:<beta}, the generic element of $Z$ has the form (\ref{eqn:k=n}). Now the proof of \cite[Theorem 5.6]{BGMN} applies to show that $G(\alpha;n,d,k)$ is irreducible.

(c) If $k>n+\frac1g(d-n)$, then, by Proposition \ref{prop:bound}, $G(\alpha;n,d,k)$ consists of a single point and (\ref{eqn:k>n}) holds. So suppose $n<k\le n+\frac1g(d-n)$. Then, by Lemma \ref{lemma:<beta}, the generic element of $Z$ has the form
\begin{equation}\label{eqn:k>n2}
0\rightarrow D^*\rightarrow V\otimes\cO\rightarrow E\rightarrow T\rightarrow0.
\end{equation}
Moreover (\ref{eqn:k>n2}) splits into two sequences
 \begin{equation}\label{eqn:vic1}
 0\rightarrow D^*\rightarrow V\otimes {\cal O}\rightarrow
 E'\rightarrow 0
 \end{equation}
and
\begin{equation}\label{eqn:vic2}
0\rightarrow E'\rightarrow E\rightarrow
 T\rightarrow 0,
 \end{equation}
where $E'$ is a vector bundle and $T$ is a torsion sheaf.

Let
$$Z':=\{(E,V)\in G(\alpha;n,d,k):(E,V) \mbox{ is generated}\}.$$
We shall prove that
$Z'$ is irreducible and that $G(\alpha;n,d,k)\setminus Z'$ is of
dimension $<\beta(n,d,k)$. Since every component of $G(\alpha;n,d,k)$ has dimension $\ge\beta(n,d,k)$, this will complete the proof.

If $(E,V)\in Z'$, then $E'=E$ in (\ref{eqn:vic1}). Dualising this sequence, we get
\begin{equation}\label{eqn:vic3}
0\rightarrow E^*\rightarrow V^*\otimes{\cO}\rightarrow D\rightarrow 0,
 \end{equation}
where $h^0(D^*)=0$ from  (\ref{eqn:vic1}) and $h^1(D)=0$ by Lemmas 
\ref{lem:h1}(a) and \ref{lem:h2}. The bundles $D$ of rank $k-n$ and degree $d$ for which $h^1(D)=h^0(D^*)=0$ form a bounded set of bundles and are therefore parametrised (not necessarily injectively) by a variety $X$ which is irreducible (or empty) 
(this follows from \cite[Theorem 2]{A}, which is essentially due to Serre; 
see also \cite[Proposition 2.6]{NR}). Let $\cD$ be the corresponding flat family over $X\times C$ and $\pi_X:X\times C\rightarrow X$ the projection. Since $H^1(D)=0$ for all $D$ in this family, $(\pi_X)_*\cD$ is a vector bundle over $X$ whose fibre over any point of $X$ corresponding to $D$ is isomorphic to $H^0(D)$. Now consider the Grassmannian bundle of subspaces $V^*$ of dimension $k$ of the fibres of $(\pi_X)_*\cD$ and the open subset $Y$ of the total space $G$ of this bundle consisting of those $V^*$ for which $(D,V^*)$ is generated. We have then an exact sequence on $Y\times C$
$$0\rightarrow\cE^*\rightarrow\cU\rightarrow(\pi\times \mbox{id}_C)^*\cD\rightarrow0,
$$
where $\cU$ is the pullback to $Y\times C$ of the universal subbundle on $G$ and $\pi:Y\rightarrow X$ is the projection. The pair $(\cE,\cU^*)$ is now a family of coherent systems on $C$ parametrised by the irreducible (or empty) variety $Y$. 
A coherent system $(E,V)$ of type $(n,d,k)$ is isomorphic to a member of this family if and only if 
$(E,V)$ is generated, $h^0(E^*)=0$ and $h^1(D)=0$, where $D^*$ is the kernel of the evaluation map $V\otimes\cO\ra E$. Since $\alpha$-stability is an open condition and $h^0(E^*)=h^1(D)=0$ for all $(E,V)\in Z'$, it follows that $Z'$ is the image of an open subset of $Y$ by some morphism and is therefore irreducible or empty.

Now let
$$Z'':=\{(E,V)\in G(\alpha;n,d,k):(E,V) \mbox{ is generically generated}\}.$$
Then $Z''$ is an open subset of $G(\alpha;n,d,k)$ consisting of those $(E,V)$ which have the form (\ref{eqn:k>n2}). By Lemma \ref{lemma:<beta}, $G(\alpha;n,d,k)\setminus Z''$ has dimension $<\beta(n,d,k)$. It is therefore sufficient to prove that $\dim(Z''\setminus Z')<\beta(n,d,k)$. In fact, if $(E,V)\in Z''\setminus Z'$ then, in the sequences (\ref{eqn:vic1}) and (\ref{eqn:vic2}),  $T$ has length $t>0$ and $\deg E'=d-t$. For fixed $t$, the extensions (\ref{eqn:vic1}) are classified by an open subset of a Quot-scheme $Q$. Tensoring (\ref{eqn:vic1}) by $D$, we see from Lemmas \ref{lem:h1}(a) and \ref{lem:h2} that $h^1(D\otimes E')=0$. It follows that the dimension of $Q$ at the point corresponding to (\ref{eqn:vic1}) is
$$h^0(D\otimes E')=k(d-t)-n(k-n)(g-1).$$
Taking account of the action of $\mbox{GL}(V)$ and (\ref{eqn:vic2}), we see that the dimension of $Z''\setminus Z'$ at $(E,V)$ is at most
$$k(d-t)-n(k-n)(g-1)+nt-(k^2-1)=\beta(n,d,k)-(k-n)t.$$
This completes the proof.

(d) For $k\le n$, it is clear from (\ref{eqn:k<n}) and (\ref{eqn:k=n}) that the Petri map is injective at the generic point of $G(\alpha;n,d,k)$. For $n<k\le n+\frac1g(d-n)$, the same follows from (\ref{eqn:k>n}) and Lemmas \ref{lem:h1}
(b) and \ref{lem:h2}. Finally, for $k>n+\frac1g(d-n)$, $G(\alpha;n,d,k)$ consists of a single point and has rank $n\le g-1$, $d=2n$ and $k=n+1$ by Proposition \ref{prop:bound}; moreover $\beta(n,2n,n+1)=0$ if and only if $n=g-1$.
\qed

\begin{corollary}\label{cor:smooth}
Suppose that $0<d\le2n$. If $G_L(n,d,k)\ne\emptyset$, then it is smooth, except possibly when $C$ is hyperelliptic and $(n,d,k)=(n,2n,n+1)$ with $n<g-1$.
\end{corollary}

\pf For $k\le n$, this is proved in \cite[Theorems 5.4 and 5.6]{BGMN}. For $n<k\le n+\frac1g(d-n)$, every element of $G_L(n,d,k)$ has the form (\ref{eqn:k>n2}) by \cite[Proposition 4.4]{BGMN}. The result follows from Lemmas 
\ref{lem:h1}(b) and \ref{lem:h2}.

Finally, suppose that $k>n+\frac1g(d-n)$. Then $G_L(n,d,k)$ consists of a single point $(E,V)$ by Proposition \ref{prop:bound} and we have an exact sequence (\ref{eqn:k>n}) with $D$ a line bundle. In the non-hyperelliptic case, $D=K$ and $E$ is a stable bundle of positive degree, so $h^1(D\otimes E)=0$. In the hyperelliptic case, $D=L^a$ and $E\cong L^{\oplus a}$ for some $a\le g-1$, where $L$ is the hyperelliptic line bundle. Under our hypotheses, this means that $a=g-1$, so $D\otimes E\cong(L^g)^{\oplus (g-1)
}$ and again $h^1(D\otimes E)=0$. It follows from (\ref{eqn:k>n}) that the Petri map of $(E,V)$ is injective; hence $G_L(g-1,2g-2,g)$ is smooth. \qed

\begin{remark}\label{rmk9}\begin{em}
If $C$ is hyperelliptic and $(n,d,k)=(n,2n,n+1)$, $n<g-1$, the Petri map cannot be injective for dimensional reasons. In this case $\beta(n,2n,n+1)<0$ and $G_L(n,2n,n+1)$ consists of the single point $D(L^n,H^0(L^n))$, but we do not know whether or not it is reduced.
\end{em}\end{remark}

\begin{corollary}\label{cor:bn}
Suppose $0<d\leq 2n$. If  $B(n,d,k)\ne\emptyset$, then it
is irreducible.
\end{corollary}

\pf Suppose $0<d\leq 2n$ and  $B(n,d,k)\ne\emptyset$. Then certainly $G_0(n,d,k)\ne\emptyset$. So, by Theorem \ref{thm:10.3}, $G_0(n,d,k)$ is irreducible. Moreover, if $g\ge3$, then $\beta(n,d,k)\le n^2(g-1)$, so \cite[Conditions 11.3]{BGMN} are satisfied and the result follows from \cite[Theorem 11.4]{BGMN}. If $g=2$ and $k>d-n$, the same argument works. If $g=2$ and $k\le d-n$, Riemann-Roch implies that  $B(n,d,k)=M(n,d)$ and is therefore irreducible. \qed

\begin{remark}\label{rmk:new}\begin{em}
Corollary \ref{cor:bn} is an improvement on results obtained in \cite{M1} 
and \cite{BGMMN}.\end{em}\end{remark}

\section{Non-emptiness}\label{section:ne}

We turn now to the question of non-emptiness of the moduli spaces. We begin by defining
$$U(n,d,k):=\{(E,V)\in G_L(n,d,k): E \mbox{ is stable}\}$$
and
$$U^s(n,d,k):=\{(E,V): (E,V) \mbox{ is $\alpha$-stable for } \alpha>0, \alpha(n-k)<d\}.$$
Note that $U(n,d,k)$ can be defined alternatively as
$$U(n,d,k):=\{(E,V): E \mbox{ is stable and } (E,V) \mbox{ is $\alpha$-stable for } \alpha>0, \alpha(n-k)<d\}$$ and in particular $U(n,d,k)\subset U^s(n,d,k)$. In the
converse direction, note that, if $(E,V)\in U^s(n,d,k)$, then $E$ is
semistable. However it is not generally true that $U(n,d,k)=U^s(n,d,k)$ and we
can have $U^s(n,d,k)\ne\emptyset$, $U(n,d,k)=\emptyset$. Our object in this section
is to determine when these sets are non-empty.

We begin with a lemma.

\begin{lemma}\label{lemma:U2}
Suppose that $0<d\le 2n$, $k>n$ and $B(n,d,k)\ne\emptyset$. Then $U(n,d,k)\ne\emptyset$.
\end{lemma}

\pf If $k>n+\frac1g(d-n)$, then, by \cite{BGN,M1,M2}, the only possibilities for $E\in B(n,d,k)$ are as follows:
\begin{itemize}
\item if $C$ is not hyperelliptic, $E\cong D(K)$;
\item if $C$ is hyperelliptic, $E\cong L$ (the hyperelliptic line bundle).
\end{itemize}
The result follows from Corollary \ref{cor:nhyp}.

Suppose now that $n<k\le n+\frac1g(d-n)$ and that $B(n,d,k)$ is non-empty. If $(E,V)$ is a coherent system with $E\in B(n,d,k)$, then $(E,V)\in G_0(n,d,k)$.  Now, by Theorem \ref{thm:10.3}(c), $G_0(n,d,k)$ is irreducible and its generic element has the form
\begin{equation}\label{eqn:k>n3}
0\rightarrow D^*\rightarrow V\otimes\cO\rightarrow E\rightarrow0,
\end{equation}
where $h^0(E^*)=0$ and $h^0(D^*)=0$; also $h^1(D)=0$ by Lemmas \ref{lem:h1}(a) and 
\ref{lem:h2}. As shown in the proof of Theorem \ref{thm:10.3}, these extensions are parametrised by an irreducible variety. By openness of stability, the generic extension (\ref{eqn:k>n3}) has $E$ stable as well as $(E,V)\in G_0(n,d,k)$. Furthermore $D$ has rank $k-n$ and its degree $d$ satisfies
$$d\ge g(k-n)+n\ge g(k-n)+\frac{d}2,$$
so $d\ge2g(k-n)$. Now any stable bundle $D$ of this rank and degree is generated by its sections and
$$h^0(D)=d-(g-1)(k-n)= d-g(k-n)-n+k\ge k,\ \ h^1(D)=0.$$
We can therefore choose a subspace $V^*$ of $H^0(D)$ of dimension $k$ which generates $D$,
giving rise to a sequence (\ref{eqn:k>n3}) for which $h^0(E^*)=h^0(D^*)=0$, $h^1(D)=0$ and $D$ is stable. 
Hence the generic extension (\ref{eqn:k>n3}) also has $D$ stable.

Finally, let us see that, for a generic extension (\ref{eqn:k>n3}), $(E,V)\in G_L(n,d,k)$ and hence $(E,V)\in U(n,d,k)$. Suppose
that $(E,V)\not\in G_L(n,d,k)$; then there
is a proper coherent subsystem $(E',V')$ such that $\frac{k'}{n'}\geq
\frac{k}{n}$. We can clearly suppose that $(E',V')$ is generically generated.  Since $(E,V)\in G_0(n,d,k)$, we have
$\frac{d'}{n'}<\frac{d}{n}$. Now let $D':=D_{V'}(E')$. We have
$(D')^* \subset D^*$, but
 $$
  \frac{d'}{k'-n'} < \frac{dn'}{n(k'-n')}  \leq \frac{dn'}{n'k-nn'}
  =\frac{d}{k-n}.
  $$
Since $\deg D^*=-d$, $\deg (D')^*\ge-d'$, this contradicts the stability of $D^*$.
\qed

We have a corresponding result for $U^s(n,d,k)$.

\begin{compl}\label{comp}
Suppose that $0<d\le2n$, $k>n$ and $G_0(n,d,k)\ne\emptyset$. Then $U^s(n,d,k)\ne\emptyset$.
\end{compl}

\pf If $k>n+\frac1g(d-n)$ and $(E,V)\in G_0(n,d,k)$, then $E$ is semistable and Proposition \ref{prop:bound} implies that $(E,V)$ is generated and $k=n+1$. The result follows from Lemma \ref{lemma:all}. If $n<k\le n+\frac1g(d-n)$, the proof of Lemma \ref{lemma:U2} still works. \qed

We are now ready to prove our main results on non-emptiness. We will state the result separately for non-hyperelliptic and hyperelliptic curves and begin with a proposition which applies in both cases.

\begin{proposition}\label{prop:non} 
Suppose $n\ge2$ and $0<d\le2n$. Then $U(n,d,k)\ne\emptyset$ if and only if one of the following three conditions applies:
\begin{itemize}
\item $0<d<2n$, $k\le n+\frac1g(d-n)$, $(n,d,k)\ne(n,n,n)$;
\item $C$ is non-hyperelliptic, $d=2n$ and either $k\le n+\frac{n}g$ or $(n,d,k)=(g-1,2g-2,g)$;
\item $C$ is hyperelliptic, $d=2n$ and $k\le n$.
\end{itemize}\end{proposition}

\pf For $k\le n$, this is proved in \cite[Theorem 3.3(v)]{BGMMN}.

For $k>n$, the stated conditions are precisely those for which $B(n,d,k)\ne\emptyset$ \cite{BGN,M1,M2} and the result follows from Lemma \ref{lemma:U2}.\qed

\begin{theorem}\label{thm:non}
Suppose that $C$ is non-hyperelliptic of genus $g\ge3$, $n\ge2$ and $0<d\le2n$. Then
$U(n,d,k)\ne\emptyset$ if and only if either 
$$k\le n+\frac1g(d-n),\ \ (n,d,k)\ne(n,n,n)$$
or 
$$(n,d,k)=(g-1,2g-2,g).$$
In all other cases,
\begin{itemize}
\item $G(\alpha;n,d,k)=\emptyset$ for all $\alpha>0$;
\item $B(n,d,k)=\emptyset$.
\end{itemize}
\end{theorem}

\pf The first part is just Proposition \ref{prop:non}. The last part follows 
from Proposition \ref{prop:bound}, except for the fact that 
$G(\alpha;n,n,n)=\emptyset$ for all $\alpha>0$, for which see 
\cite[Remark 5.7]{BGMN}.\qed

\begin{theorem}\label{thm:hyp}
Suppose that $C$ is hyperelliptic, $n\ge2$ and $0<d\le2n$. Then
\begin{itemize}
\item[(a)] $U(n,d,k)\ne\emptyset$ if and only if either 
$$0<d<2n,\ \ k\le n+\frac1g(d-n),\ \ (n,d,k)\ne(n,n,n)$$
or $d=2n, k\le n$;
\item[(b)] if $k>n$, then
\begin{itemize}
\item $U(n,2n,k)=\emptyset$;
\item $U^s(n,2n,k)\ne\emptyset$ if and only if either $k\le n+\frac{n}g$ or $k=n+1$ and $2\le n\le g-1$.
\end{itemize}
\end{itemize}
In all other cases,
\begin{itemize}
\item $G(\alpha;n,d,k)=\emptyset$ for all $\alpha>0$;
\item $B(n,d,k)=\emptyset$.
\end{itemize}
\end{theorem}

We already have enough information to prove this except for showing that $U^s(n,2n,k)\ne\emptyset$ when $n<k\le n+\frac{n}g$. This will be done in the next section.

\begin{remark}\begin{em} The case $n=1$ has been explicitly excluded from
these statements as the results need modification. 
In this case the $\alpha$-stability condition is redundant and the triples 
for which $0<d\le2$ and $U(1,d,k)\ne\emptyset$ are $(1,1,1)$, $(1,2,1)$
and, for $C$ hyperelliptic, $(1,2,2)$.\end{em}\end{remark} 

\section{Proof of Theorem \ref{thm:hyp}}\label{section:hyp}

In this section we suppose that $C$ is hyperelliptic and $L$ is the hyperelliptic line bundle. We assume that $n\ge2$ and investigate by a sequence of propositions the case 
\begin{equation}\label{eqn:hyp1}
d=2n,\ \ n<k\le n+\frac{n}g.
\end{equation}
\begin{proposition}\label{prop:hyp1}
Suppose $C$ is hyperelliptic. Then $U^s(n,2n,n+1)\ne\emptyset$.
\end{proposition}
\pf Let $E=L^{\oplus n}$. Then $E$ is generated and we can choose a subspace $V$ of $H^0(E)$ of dimension $n+1$ such that $(E,V)$ is generated. The result follows from Lemma \ref{lemma:all}.\qed

\begin{remark}\label{rmk:hyp1}
\begin{em} Proposition \ref{prop:hyp1} applies even when $n+1>n+\frac{n}g$, in which case it has already been proved in Corollary \ref{cor:hyp}\end{em}.
\end{remark}

Now suppose that $k\ge n+2$ and write $k=n+r$, so that (\ref{eqn:hyp1}) becomes
$$d=2n,\ \ 2\le r\le\frac{n}g.$$

\begin{proposition}\label{prop:hyp2}
Suppose $C$ is hyperelliptic and $2\le r\le\frac{n-2}g$. Then $U^s(n,2n,n+r)\ne\emptyset$.
\end{proposition}
\pf We consider extensions
\begin{equation}\label{eqn:hyp2}
0\lra(E_1,V_1)\lra(E,V)\lra(E_2,V_2)\lra0,
\end{equation}
where $(E_1,V_1)$ has type
$$(n_1,d_1,k_1)=(n-1,2n-3,n+r-1)$$
and $(E_2,V_2)$ has type $(1,3,1)$. Certainly $(E_2,V_2)\in U(1,3,1)$. On the other hand $d_1<2n_1$ and
\begin{equation}\label{eqn:e1}
k_1=n_1+r\le n_1+\frac{n-2}g=n_1+\frac1g(d_1-n_1).
\end{equation}
So, by Proposition \ref{prop:non}, we can choose $(E_1,V_1)\in U(n_1,d_1,k_1)$. 
To show that there exist non-trivial extensions (\ref{eqn:hyp2}), it is sufficient to prove that $C_{21}>0$. In fact, by (\ref{eqn:basic4}),
\begin{eqnarray}\label{eqn:extra}
C_{21}&=& n_1(n_2-k_2)(g-1)+(k_2-n_2)d_1+d_2n_1-k_1k_2\nonumber\\
&=&3(n-1)-(n+r-1)=2n-2-r\ge2n-2-\frac{n-2}g>0.
\end{eqnarray}

Suppose now that (\ref{eqn:hyp2}) is non-trivial. If $\alpha_c=\frac{n}r$, then
$$\mu_{\alpha_c}(E_1,V_1)=\frac{2n-3}{n-1}+\frac{n}r\cdot\frac{n+r-1}{n-1}=3+\frac{n}r=\mu_{\alpha_c}(E_2,V_2).$$
Since $\mu_{\alpha^-_c}(E_1,V_1)<\mu_{\alpha^-_c}(E_2,V_2)$, it follows 
from Proposition \ref{prop:basic} that $(E,V)$ is $\alpha^-_c$-stable.

Now consider the extension of bundles
\begin{equation}\label{eqn:hyp3}
0\lra E_1\lra E\lra E_2\lra0
\end{equation}
underlying (\ref{eqn:hyp2}) and suppose first that this extension is non-trivial. If $F$ is a subbundle of $E$ which contradicts semistability, then certainly $F\not\subset E_1$. Moreover, if $F\ra E_2$ is not surjective, then we have an extension
$$0\lra F_1\lra F\lra F_2\lra0$$
with $\deg F_2\le2$ and $F_1\subset E_1$, so $\mu(F)<2$. It follows that, to contradict semistability of $E$, we must have
$$0\lra F_1\lra F\lra E_2\lra0.$$
Moreover $F_1\ne0$ since (\ref{eqn:hyp3}) does not split. Since $E_1$ is stable, $\mu(F_1)<\mu(E_1)<2$, so
$$\deg F_1\le2\,\mbox{rk}F_1-1,\ \ \deg F=\deg F_1+3\le2\,\mbox{rk}F.$$
So $E$ is semistable.

To complete the proof in this case, it is sufficient by Complement \ref{comp} to show that $(E,V)\in U_0(n,2n,n+r)$. If this is not the case, there exists a proper coherent subsystem $(F,W)$ of type $(n_F,d_F,k_F)$ of $(E,V)$ with $d_F=2n_F$ and $\frac{k_F}{n_F}\ge\frac{n+r}n$. But in this case $(E,V)$ cannot be $\alpha$-stable for any $\alpha>0$, contradicting the fact that $(E,V)$ is $\alpha^-_c$-stable. 

It remains to prove that there exist extensions (\ref{eqn:hyp2}) such that (\ref{eqn:hyp3}) does not split. Now, by \cite[Corollaire 1.6]{He} (see also \cite[equation (7)]{BGMN}), we have an exact sequence
$$\Hom(V_2,H^0(E_1)/V_1)\lra\Ext^1((E_2,V_2),(E_1,V_1))\lra\Ext^1(E_2,E_1).$$
It is therefore sufficient to prove that

$$\dim\Ext^1((E_2,V_2),(E_1,V_1))>\dim\Hom(V_2,H^0(E_1)/V_1).$$
Now, by (\ref{eqn:extra}),
$$\dim\Ext^1((E_2,V_2),(E_1,V_1))\ge C_{21}=2n-2-r,$$
while
$$\dim\Hom(V_2,H^0(E_1)/V_1)=h^0(E_1)-(n+r-1).$$
By \cite[Chapitre 2, Th\'eor\`eme A.1]{M1}, we have
$$h^0(E_1)\le n_1+\frac1g(d_1-n_1)=n-1+\frac{n-2}g,$$
so 
$$\dim\Hom(V_2,H^0(E_1)/V_1)\le\frac{n-2}g-r<C_{21}\le\dim\Ext^1((E_2,V_2),(E_1,V_1)).$$
\qed

\begin{remark}\label{rmk:hyp2}\begin{em}
It is perhaps of interest to note that the coherent systems $(E,V)$ constructed in this proof are not themselves in $U^s(n,2n,n+r)$. We need to use Complement \ref{comp} to prove the proposition. Moreover the hypothesis $r\le\frac{n-2}g$ is used in an essential way (see (\ref{eqn:e1})) and the method of proof does not work without it; in fact, without the hypothesis, there are no flips.
\end{em}\end{remark}

It remains to consider the cases $r=\frac{n-1}g$ and $r=\frac{n}g$. In other words, we have two cases
$$
n=gr+1,\ \ r\ge2
$$
and
$$
n=gr,\ \ r\ge2.
$$

\begin{proposition}\label{prop:hyp3}
Suppose $C$ is hyperelliptic and $r\ge2$. Then $U^s(gr+1,2gr+2,gr+r+1)\ne\emptyset$.
\end{proposition}

\pf We consider extensions
\begin{equation}\label{eqn:hyp4} 
0\lra(E_1,V_1)\lra(E,V)\lra(E_2,V_2)\lra0,
\end{equation}
where $(E_2,V_2)\cong D(L^{g-1},H^0(L^{g-1}))$ and 
$$(E_1,V_1)\in U^s(g(r-1)+2,2g(r-1)+4,g(r-1)+r+1),$$
which is non-empty by Propositions \ref{prop:hyp1} and \ref{prop:hyp2}. By 
Theorem \ref{thm:10.3}(c), we can suppose further that $(E_1,V_1)$ is 
generated.
Note also that $(E_2,V_2)$ is generated and has the form $(L^{\oplus(g-1)},V_2)$ with $\dim V_2=g$, 
and belongs to $U^s(g-1,2g-2,g)$ by Corollary \ref{cor:hyp}.

We show first that there exists a non-trivial extension (\ref{eqn:hyp4}). 
In fact, by (\ref{eqn:basic4}),
\begin{eqnarray*}
C_{21}&=& n_1(n_2-k_2)(g-1)+(k_2-n_2)d_1+d_2n_1-k_1k_2\\
&=&-n_1(g-1)+2n_1g-k_1g\\&=&n_1+g(n_1-k_1)=g(r-1)+2-g(r-1)=2.
\end{eqnarray*}

From now on we suppose that (\ref{eqn:hyp4}) is non-trivial. Let $(E',V')$ 
be a coherent subsystem of $(E,V)$ of type $(n',d',n'+r')$ which 
contradicts $0^+$-stability. Then certainly $E'$ is semistable of slope $2$, so $d'=2n'$, and
\begin{equation}\label{eqn:hyp5}
\frac{r'}{n'}\ge\frac{r}n=\frac{r}{gr+1}.
\end{equation}
From (\ref{eqn:hyp4}), we have an extension
\begin{equation}\label{eqn:hyp6}
0\lra(E_1',V_1')\lra(E',V')\lra(E_2',V_2')\lra0.
\end{equation}
Since $E_1$, $E_2$ are semistable of slope $2$, so are $E_1'$ and $E_2'$. 
For $i=1,2$, denote the type of $(E_i',V_i')$ by $(n_i',2n_i',n_i'+r_i')$. 
Note that $n_1'\ne0$, for otherwise (\ref{eqn:hyp6}) would require either $r'\le0$, contradicting  
(\ref{eqn:hyp5}), or $(E',V')=(E_2,V_2)$, in which case it 
would split the sequence (\ref{eqn:hyp4}). Since $(E_1,V_1)$ is $0^+$-stable,
we have
\begin{equation}\label{eqn:hyp7}
\frac{r_1'}{n_1'}\le\frac{r_1}{n_1}=\frac{r-1}{g(r-1)+2}.
\end{equation}
If $(E_2',V_2')\ne(E_2,V_2)$, then $r_2'\le0$; this, together with 
(\ref{eqn:hyp7}), contradicts (\ref{eqn:hyp5}). Hence (\ref{eqn:hyp6}) becomes
\begin{equation}\label{eqn:600}
0\lra(E_1',V_1')\lra(E',V')\lra(E_2,V_2)\lra0,
\end{equation}
from which it follows that 
$$n_1'=n'-g+1,\ \ r_1'=r'-1.$$
A simple calculation shows that equations (\ref{eqn:hyp5}) and 
(\ref{eqn:hyp7}) can be written as
\begin{equation}\label{eqn:hyp8}
n_1'\le r_1'g+1+\frac{r'}r
\end{equation}
and
\begin{equation}\label{eqn:hyp9}
n_1'\ge r_1'g+\frac{2r_1'}{r-1}.
\end{equation}
Since $r'<r$, (\ref{eqn:hyp8}) implies that $n_1'\le r_1'g+1$. By equation 
(\ref{eqn:hyp9}), this is only possible when
\begin{equation}\label{eqn:hyp9'}
n_1'=r_1'g+1,\ \ r_1'\le\frac{r-1}2.
\end{equation}  

Note that $(E_1',V_1')$ is a coherent subsystem of $(E_1,V_1)$. We must
have $V_1'=V_1\cap H^0(E_1')$, otherwise we could replace $V_1'$ by a 
subspace of $H^0(E_1')$ of greater dimension, which would contradict 
(\ref{eqn:hyp9}). Thus we have an extension
\begin{equation}\label{eqn:hyp10}
0\lra(E_1',V_1')\lra(E_1,V_1)\lra(F,W)\lra0.
\end{equation}
We now count parameters to show that the $(E_1,V_1)$ occurring in an extension
(\ref{eqn:hyp10}) are not generic.

We begin with two lemmas.

\begin{lemma}\label{lem:a1}
$(E_1',V_1')$ is generically generated and  $0^+$-stable.
\end{lemma}

\pf If $r_1'=0$, then, by (\ref{eqn:hyp9'}), $(E_1',V_1')$ has type 
$(1,2,1)$ and the result is immediate. So suppose $r_1'\ge1$.

If $(E_1',V_1')$ is not generically generated, it possesses a coherent 
subsystem $(E_1'',V_1')$ with $\rk E_1''\le n_1'-1=r_1'g$. 
By $\alpha$-stability of $(E_1,V_1)$ for large $\alpha$, 
this implies
$$\frac{r_1'g+1+r_1'}{r_1'g}\le\frac{g(r-1)+r+1}{g(r-1)+2},$$
which is evidently false. On the other hand, if $(E_1',V_1')$ is not 
$0^+$-stable, there exists a proper coherent subsystem $(E_1'',V_1'')$ of 
type $(n_1'',2n_1'',n_1''+r_1'')$ such that
$$\frac{r_1''}{n_1''}\ge\frac{r_1'}{n_1'}=\frac{r_1'}{r_1'g+1},$$
i.~ e.
\begin{equation}\label{eqn:hyp11}
n_1''\le r_1''\left(g+\frac1{r_1'}\right).
\end{equation}
By $\alpha$-stability of $(E_1,V_1)$, we have also
$$\frac{r_1''}{n_1''}<\frac{r-1}{g(r-1)+2},$$
i.~ e.
\begin{equation}\label{eqn:hyp12}
n_1''> r_1''\left(g+\frac2{r-1}\right).
\end{equation}
Now $n_1''<n_1'$, so
$$r_1''\left(g+\frac2{r-1}\right)< n_1'-1=r_1'g.$$
Hence $r_1''<r_1'$ and (\ref{eqn:hyp11}) and (\ref{eqn:hyp12}) give a 
contradiction.\qed

\begin{lemma}\label{lem:a2} $(F,W)$ is generated and
$\Hom(D(K,H^0(K)),(F,W))=0$.
\end{lemma} 

\pf Since $(E_1,V_1)$ is generated, so is $(F,W)$. Now suppose that $\phi:D(K,H^0(K))\ra(F,W)$ is a non-zero homomorphism. 
Since 
$D(K,H^0(K))$ is $\alpha$-stable for $\alpha>0$, and $F$ is semistable of 
slope $2$, the image of $\phi$ is a coherent subsystem $(F',W')$ of type 
$(n_{F'},2n_{F'},n_{F'}+r_{F'})$ of $(F,W)$ with
\begin{equation}\label{eqn:hyp13}
n_{F'}\le g-1,\ \ r_{F'}\ge1.
\end{equation}
The pullback of $(F',W')$ to $(E_1,V_1)$ in (\ref{eqn:hyp10}) has type
$$(n_{F'}+n_1',2(n_{F'}+n_1'),n_{F'}+r_{F'}+n_1'+r_1').$$
Now $\alpha$-stability of $(E_1,V_1)$ gives
$$\frac{r_{F'}+r_1'}{n_{F'}+n_1'}\le\frac{r-1}{g(r-1)+2},$$
i.~e.
$$n_{F'}+n_1'\ge(r_{F'}+r_1')\left(g+\frac2{r-1}\right).$$
Since $n_1'=r_1'g+1$ by (\ref{eqn:hyp9'}), this is equivalent to
$$n_{F'}+1\ge r_{F'}\left(g+\frac2{r-1}\right)+\frac{2r_1'}{r-1}.$$
This contradicts (\ref{eqn:hyp13}).\qed

For our parameter count, we now establish three claims.

\begin{claim}\label{claim1}
For fixed $(E_1',V_1')$, $(F,W)$, the non-trivial extensions (\ref{eqn:hyp10})
for which $(E_1,V_1)$ is generated and $\alpha$-stable for some $\alpha$ 
depend on at most
$$C_{21}^{(\ref{eqn:hyp10})}-\dim\Aut(F,W)$$
parameters, where $C_{21}^{(\ref{eqn:hyp10})}$ denotes the value of $C_{21}$ 
for the extensions (\ref{eqn:hyp10}).
\end{claim}

\pf By (\ref{eqn:basic3}), we have
$$\dim\Ext^1((F,W),(E_1',V_1'))=C_{21}^{(\ref{eqn:hyp10})}+\dim\HH^0_{21}
+\dim\HH^2_{21},$$
where
$$\HH^0_{21}=\Hom((F,W),(E_1',V_1')),\ \ \dim\HH^2_{21}=\Ext^2((F,W),
(E_1',V_1')).$$
Now $\HH^0_{21}=0$ since otherwise (\ref{eqn:hyp10}) would give a 
contradiction to the $\alpha$-stability of $(E_1,V_1)$. On the other hand, 
by (\ref{eqn:basic7}) and Serre duality,
$$\HH^2_{21}\cong H^1(E_1'\otimes N_2^*),$$
where $N_2$ is defined by an exact sequence
$$0\lra N_2\lra W\otimes \cO\lra F\lra0.$$
(Note that $(F,W)$ is generated by Lemma \ref{lem:a2}.) By Lemmas \ref{lem:a2}
and \ref{lem:h1}(a), we have $h^1(N_2^*)=0$. By Lemma \ref{lem:a1}, we have an 
exact sequence
$$0\lra N_1\lra V_1'\otimes\cO\lra E_1'\lra T_1\lra0,$$
where $T_1$ is a torsion sheaf. Hence $h^1(E_1'\otimes N_2^*)=0$. Finally, 
since we are assuming $(E_1,V_1)$ is $\alpha$-stable for some $\alpha$, the 
action of $\Aut(F,W)$ on the extensions (\ref{eqn:hyp10}) is free. The result 
follows.\qed

\begin{claim}\label{claim2}
$(E_1',V_1')$ depends on at most
$$\beta(n_1',2n_1',k_1')$$
parameters.
\end{claim}

\pf Since $(E_1',V_1')$ is $0^+$-stable by Lemma \ref{lem:a1}, it is 
sufficient to show that the Petri map is 
injective. This follows from Lemmas \ref{lem:h1}(b) and \ref{lem:h2} and 
(\ref{eqn:hyp9'}).\qed 

\begin{claim}\label{claim3}
$(F,W)$ depends on at most
$$\beta(n_1-n_1',2n_1-2n_1',k_1-k_1')+\dim\Aut(F,W)-1$$
parameters.
\end{claim}

\pf By Lemmas \ref{lem:a2} and \ref{lem:h1}(b), the Petri map of $(F,W)$ is injective and so, 
by Proposition \ref{prop:extra}, $\Ext^2((E,V),(E,V))=0$. It follows from (\ref{eqn:basic8}) that any $(F,W)$ 
has a local deformation space of dimension
 $$\dim\Ext^1((F,W),(F,W))=\beta(n_1-n_1',2n_1-2n_1',k_1-k_1')+
\dim\Aut(F,W)-1.$$
The result follows.\qed

\begin{em}Completion of proof of Proposition \ref{prop:hyp3}.\end{em}
By the above claims, the $(E_1,V_1)$ for which 
$(E_1',V_1')$ exists satisfying (\ref{eqn:hyp8}) and (\ref{eqn:hyp9}) depend 
on at most
$$\beta(n'_1,2n'_1,k'_1)+\beta(n_1-n'_1,2n_1-2n'_1,k_1-k'_1)+C^{(\ref{eqn:hyp10})}_{21}-1$$
parameters. By (\ref{eqn:basic6}), this number is equal to
$$\beta(n_1,2n_1,k_1)-C_{12}^{(\ref{eqn:hyp10})}.$$
If $C_{12}^{(\ref{eqn:hyp10})}>0$, it follows that the general 
$(E_1,V_1)$ contains no $(E_1',V_1')$ with these properties. Hence $(E,V)$ 
contains no subsystem $(E',V')$ contradicting $0^+$-stability. The result 
then follows from Complement \ref{comp}.

It remains to calculate $C_{12}^{(\ref{eqn:hyp10})}$. In fact, by 
(\ref{eqn:basic5}),
\begin{eqnarray*}
C_{12}^{(\ref{eqn:hyp10})}&=&(n_1-n_1')(n_1'-k_1')(g-1)+2(n_1-n_1')(k_1'-n_1')
+2n_1'(n_1-n_1')-k_1'(k_1-k_1')\\
&=&(n_1-n_1')[(n_1'-k_1')(g-1)+2k_1']-k_1'(k_1-k_1')\\
&=&(g(r-1-r_1')+1)(-r_1')(g-1)+k_1'(2n_1-2n_1'-k_1+k_1')\\
&=&-(r-1-r_1')(g-1)r_1'g-r_1'(g-1)+(r_1'(g+1)+1)((g-1)(r-1-r_1')+1)\\
&=&(r-1-r_1')(g-1)(r_1'+1)+2r_1'+1>0.
\end{eqnarray*}\qed

\begin{proposition}\label{prop:hyp4}
Suppose $C$ is hyperelliptic and $r\ge2$. Then $U^s(gr,2gr,gr+r)\ne
\emptyset$.
\end{proposition}

\pf The proof is similar to that of Proposition \ref{prop:hyp3}; we outline below the necessary changes.

We consider sequences (\ref{eqn:hyp4}), where $(E_2,V_2)\cong D(K,H^0(K))$
as before, and now
$$(E_1,V_1)\in U^s(g(r-1)+1,2g(r-1)+2,g(r-1)+r).$$
This space is non-empty by Propositions \ref{prop:hyp1} and \ref{prop:hyp3}. 
We have
$$C_{21}=n_1+g(n_1-k_1)=1,$$
so non-trivial extensions (\ref{eqn:hyp4}) exist. After modifying 
(\ref{eqn:hyp5}) and (\ref{eqn:hyp7}), we proceed to (\ref{eqn:hyp8}) and 
(\ref{eqn:hyp9}), which become
$$n_1'\le r_1'g+1,\ \ n_1'\ge r_1'g+\frac{r_1'}{r-1}.$$
So again $n_1'=r_1'g+1$, where now $r_1'\le r-1$. In fact, since $(E_1',V_1')\ne(E,V)$ in (\ref{eqn:600}) and $(E_1,V_1)$ is generated, $r_1'<r_1=r-1$. 

The proofs of Lemmas \ref{lem:a1} and \ref{lem:a2} are the same as before, 
replacing $g(r-1)+2$ by $g(r-1)+1$ with consequential changes which don't 
affect the argument. The only remaining thing to be checked is that 
$C_{12}^{(\ref{eqn:hyp10})}>0$. In fact
\begin{eqnarray*}
C_{12}^{(\ref{eqn:hyp10})}&=&(n_1-n_1')[(n_1'-k_1')(g-1)+2k_1']-k_1'
(k_1-k_1')\\
&=&g(r-1-r_1')(-r_1')(g-1)+k_1'(2n_1-2n_1'-k_1+k_1')\\
&=&-(r-1-r_1')(g-1)r_1'g+(r_1'(g+1)+1)(g-1)(r-1-r_1')\\
&=&(r-1-r_1')(g-1)(r_1'+1)>0.
\end{eqnarray*}\qed

{\em Proof of Theorem \ref{thm:hyp}.}
(a) is just Proposition \ref{prop:non}. (b) follows from Propositions \ref{prop:hyp1}, \ref{prop:hyp2}, \ref{prop:hyp3} and \ref{prop:hyp4} in the case $k\le n+\frac{n}g$ and from Corollary \ref{cor:hyp} if $k>n+\frac{n}g$. The last part follows from Proposition \ref{prop:bound}.\qed

\section{An example with $d>2n$}\label{section:ex}

We have seen that, when $d<2n$ and $k>n$ and $\alpha$-stable
coherent systems exist for some $\alpha$ (i.~e. when
(\ref{eq:kbound}) holds), then there exist coherent systems
$(E,V)$ such that $E$ is stable and $(E,V)$ is $\alpha$-stable for
all $\alpha>0$. The same applies when $d=2n$ if $C$ is not
hyperelliptic. If $C$ is hyperelliptic of genus $g\ge3$ and
$a\ge2$, the coherent systems $(L^{\oplus a},W)$ of type
$(a,2a,a+1)$ are $\alpha$-stable for all $\alpha$ by Corollary
\ref{cor:hyp}, but $L^{\oplus a}$ is only semistable. Moreover,
when $d\le 2n$, there is no case in which there exist semistable
bundles but $\alpha$-stable coherent systems do not exist for
large $\alpha$. The object of this section is to construct such
examples, necessarily with $d>2n$.

\begin{lemma}\label{lemma:notlarge}
Suppose $(E,V)$ is a coherent system of type $(n,d,k)$ with $h^0(E^*)=0$ and 
\be\label{eq:nonss} n+\frac{1}{g}(d-n)<k<\frac{ng}{g-1}.
\end{equation}
Then $(E,V)$ is not $\alpha$-semistable for large $\alpha$.
\end{lemma}

\pf Suppose (\ref{eq:nonss}) is satisfied. By Corollary
\ref{cor:stable}, there exists a non-zero homomorphism
$$
D(K,H^0(K))\longrightarrow (E,V).
$$
By Corollaries \ref{cor:nhyp} and \ref{cor:hyp}, $D(K,H^0(K))$ is $\alpha$-stable for all $\alpha>0$. Moreover
$$
\mu_{\alpha}(D(K,H^0(K)))=2+\alpha\frac{g}{g-1}>\frac{d}{n}+\alpha\frac{k}n
$$
for sufficiently large $\alpha$ by (\ref{eq:nonss}). This
contradicts the $\alpha$-semistability of $(E,V)$ for large
$\alpha$. \qed

\begin{proposition}\label{prop:ex}
Suppose that $C$ is not hyperelliptic and $3\le r\le g-1$. Then
there exists a coherent system $(E,V)$ of type \be\label{eq:type}
(rg-r+1,2rg-2r+3,rg+1)
\end{equation}
with $E$ stable. Moreover (\ref{eq:nonss}) is satisfied.
\end{proposition}

\pf It is clear that (\ref{eq:nonss}) follows from (\ref{eq:type})
and the assumption $r\ge3$.

Consider the extensions \be\label{eq:ext} 0\longrightarrow
D(K)^{\oplus r}\oplus \cO(p_1,p_2)\longrightarrow
E\longrightarrow {\cal O}_q\longrightarrow 0,
\end{equation}
where $p_1,p_2,q\in C$. We take $V$ to be the image of
$H^0(D(K)^{\oplus r}\oplus \cO(p_1,p_2))$ in $H^0(E)$. The fact
that $(E,V)$ is of type (\ref{eq:type}) is clear. Since $D(K)$ and
$\cO(p_1,p_2)$ are both stable of slope $2$, $E$ fails to be
stable only if it admits either $D(K)$ or $\cO(p_1,p_2)$ as a
quotient. Equivalently $E$ fails to be stable only if some factor
of $D(K)$ or $\cO(p_1,p_2)$ splits off (\ref{eq:ext}). Now the
extensions (\ref{eq:ext}) are classified by elements
\begin{eqnarray*}
e=(e_1,\ldots,e_r,e_{r+1})&\in&\Ext^1(\cO_q,D(K)^{\oplus r}\oplus \cO(p_1,p_2))\\
&&=\Ext^1(\cO_q,D(K))^{\oplus r}\oplus\Ext^1(\cO_q,\cO(p_1,p_2)).
\end{eqnarray*}
Since $D(K)$ and $\cO(p_1,p_2)$ are stable and non-isomorphic, it
follows that $E$ is stable if $e_1,\ldots,e_r$ are linearly
independent and $e_{r+1}\ne0$. Since
$$
\dim \Ext^1({\cal O}_q,D(K))= g-1,
$$
it is possible to choose such $e_1,\ldots e_r,e_{r+1}$ whenever
$r\le g-1$. This completes the proof. \qed

\section*{Acknowledgments}
The fourth author wishes to acknowledge Universidad Complutense de Madrid and the Institute for Advanced Study, Princeton for their hospitality and for providing excellent working conditions. The fifth author would like to thank CIMAT, Guanajuato, Mexico and California State University Channel Islands, where parts of this research were carried out.

\end{document}